\documentclass[11pt]{article}
\usepackage{latexsym,amsthm,amsmath,amssymb,graphicx,enumitem}

\usepackage[utf8]{inputenc}

\newtheorem{theorem}{Theorem}[section]
\newtheorem{lemma}[theorem]{Lemma}
\newtheorem{corollary}[theorem]{Corollary}

\newcommand{\parag}[1]{\vspace{2mm}

\noindent{\bf #1} }

\newcommand{\RR}{\ensuremath{\mathbb R}}
\newcommand{\pts}{\mathcal P}
\newcommand{\lines}{\mathcal L}
\newcommand{\circs}{\mathcal C}
\newcommand{\curves}{\Gamma}
\def\eps{{\varepsilon}}
\newcommand{\bis}{\mathcal B}
\newcommand{\crs}{\mathrm{cr}}

\title{Distinct Distances with $\ell_p$ Metrics}
\author{Polly Matthews Jr.\footnote{This work was done as part of the \emph{Polymath REU} program. The authors: Moaaz AlQady (The American University in Cairo), Riley Chabot (Princeton University), William Dudarov (Carleton College), Linus Ge (University of Rochester), Mandar Juvekar (University of Rochester), Srikanth Kundeti (Rutgers University), Neloy Kundu (Lafayette College), Kevin Lu (Georgia Institute of Technology), Yago Moreno (University of Bristol), Sibo Peng (North Carolina State University), Samuel Speas (University of California, Berkeley), Julia Starzycka (University of Illinois at Chicago), Henry Steinthal (Lafayette College), and Anastasiia Vitko (Wesleyan University). }}

\date{}

\begin{document}

\maketitle

\begin{abstract}
We study Erd\H os's distinct distances problem under $\ell_p$ metrics with integer $p$.
We improve the current best bound for this problem from $\Omega(n^{4/5})$ to  $\Omega(n^{6/7-\eps})$, for any $\eps>0$. We also characterize the sets that span an asymptotically minimal number of distinct distances under the $\ell_1$ and $\ell_\infty$ metrics.
\end{abstract}

\section{Introduction}

Erd\H os proposed the \emph{distinct distances problem} in his seminal 1946 paper \cite{erd46}.
For a finite point set $\pts\subset \RR^2$, let $D(\pts)$ be the set of distances spanned by pairs of points from $\pts$.
The distinct distances problem asks for the asymptotic value of $\min_{|\pts|=n} |D(\pts)|$.
In other words, the problem asks for the minimum number of distances that could be spanned by $n$ points in $\RR^2$.

Erd\H os \cite{erd46} presented a set of $n$ points that spans $\Theta(n/\sqrt{\log n})$ distinct distances. He also conjectured that no set spans an asymptotically smaller number of distances.
Over the decades, a large number of works have been dedicated to this problem. 
Recently, Guth and Katz \cite{GK15} almost completely settled Erd\H os's conjecture, proving that $D(\pts)=\Omega(n/\log n)$ for any set $\pts\subset\RR^2$ of $n$ points. 

Over the decades, the distinct distances problem evolved into a large family of problems, mostly posed by Erd\H os. Even after the Guth--Katz breakthrough, most of the main distinct distances problems remain wide open (for example, distinct distances in $\RR^d$ and the structural problem in $\RR^2$).   
A deep distinct distances theory is being developed and an entire book is dedicated to it \cite{GIS11}.
For a survey of distinct distances problems in $\RR^d$, see \cite{Sheffer14}.
The problem has also been studied in complex spaces \cite{SheZa20}, spaces over finite fields \cite{BKT04,MPPRS03}, and more. 

Instead of changing the underlying field, one may change the distance metric. A result of Matou\v sek \cite{Mat11} implies that, for ``most'' metrics in $\RR^2$, every set of $n$ points in $\RR^2$ determines $\Omega(n/(\log n \cdot \log \log n))$ distinct distances.
Since this is not relevant for the current paper, we do not clarify the exact meaning of ``most''.

In the current work we focus on $\ell_p$ metrics. In other words, we consider the metrics induced by the $\ell_p$ norms.
In particular, we focus on integer $p$.
For an integer $p\ge 1$, the $\ell_p$ distance between the points $(a_x,a_y)$ and $(b_x,b_y)$ is 
\[ \left(|a_x-b_x|^p+|a_y-b_y|^p\right)^{1/p}. \]

\parag{New bounds for $\ell_p$ metrics}Distinct distances under various metrics has been the topic of the dissertation of Julia Garibaldi \cite{Garibaldi04} (under the supervision of Terence Tao).
The following result is obtained in Garibarldi's dissertation.
Given a point $u$ and a point set $\pts$, we denote by $D(u,\pts)$ the number of distances between $u$ and the points of $\pts$.

\begin{theorem} \label{th:GariLp}
Let $\pts$ be a set of $n$ points in $\RR^2$ and let $p> 1$. Then, under the $\ell_p$ metric, there exists a point $u\in \pts$ such that 
\[ D(u,\pts) = \Omega\left(n^{4/5}\right).\]
\end{theorem}

When $u\in \pts$, we have that $D(\pts)\ge D(u,\pts)$. Thus, the statement of Theorem \ref{th:GariLp} is stronger than a bound for the minimum number of distinct distances. 

We derive the following stronger bound for the case of integer $p$.

\begin{theorem}  \label{th:SzekelyImprov}
Let $\pts$ be a set of $n$ points in $\RR^2$ and let $p> 1$ be an integer. Then, under the $\ell_p$ metric and for any $\eps>0$, there exists a point $u\in \pts$ such that 
\[ D(u,\pts) = \Omega_{p,\eps}\left(n^{6/7-\eps}\right).\]
\end{theorem}

The proofs of Theorems \ref{th:GariLp} and \ref{th:SzekelyImprov} adapt Szekely's proof for distinct distances under the Euclidean metric \cite{Szek97}. The original proof led to the bound $D(u,\pts) = \Omega\left(n^{4/5}\right)$, and Garibaldi extended this to all $\ell_p$ metrics (and to other metrics). Surprisingly, we use a similar approach to obtain the stronger bound $D(u,\pts) = \Omega_{p,\eps}\left(n^{6/7-\eps}\right)$. We do not know how to use our approach with the Euclidean metric.

Theorem \ref{th:SzekelyImprov} was already known for the case of $p=2$ (the Euclidean distance). The Guth-Katz result does not extend to the stronger formulation involving a single point $u$.
However, a bound of $D(u,\pts) = \Omega\left(n^{6/7}\right)$ was derived by Solymosi and T\'oth \cite{SoTo01} (this was further improved by Tardos \cite{Tardos03}). While both Theorem \ref{th:SzekelyImprov} and the Solymosi-T\'oth work achieve the exponent 6/7, the two proofs are rather different. 

\parag{Structure for the cases of $\ell_1$ and $\ell_\infty$.}
Under the $\ell_\infty$ metric, the distance between the points $(a_x,a_y)$ and $(b_x,b_y)$ is 
\[ \max\{|a_x-b_x|,|a_y-b_y|\}. \]
The $\ell_1$ and $\ell_\infty$ metrics are degenerate in the sense that  their unit circles are not strictly convex (a \emph{unit circle} is the set of points at distance of one from a fixed point. These curves are not circles in the standard definition under any non-Euclidean metric.)
For distinct distances in $\RR^2$, the two metrics are equivalent. See Section \ref{sec:Structure} for more details. 
Without loss of generality, we now consider the $\ell_\infty$ metric.
The set $\pts = \{1,2,3,\ldots,\sqrt{n}\}^2$ satisfies $D(\pts)=\sqrt{n}-1$ and it is not difficult to show that no set of $n$ points satisfies $D(\pts)<\sqrt{n}-1$. (See \cite{BELMPRT19} for a generalization to $\RR^d$ for $\ell_1$.) 

One of the main open distinct distances problems is characterizing the sets $\pts\subset \RR^2$ that satisfy $D(\pts) =O(n/\sqrt{\log n})$ (under the Euclidean metric). 
This is an unusually difficult problem, with hardly anything known about it.

We study the structural problem for the $\ell_\infty$ norm (so also for the $\ell_1$ norm). In particular, we characterize the sets of $n$ points $\pts\subset \RR^2$ that satisfy $D(\pts)=\Theta(\sqrt{n})$.
For a definition of generalized arithmetic progressions, see Section \ref{sec:Structure}.

\begin{theorem} \label{th:Structure}
Let $\pts$ be a set of $n$ points such that $D(\pts)=\Theta(\sqrt{n})$ under the $\ell_\infty$ metric. Then: \\[2mm]
(a) There exists a set $\lines$ of $\Theta(\sqrt{n})$ lines such that $\pts \subset \bigcup_{\ell\in \lines} \ell$. Either all lines of $\lines$ are horizontal, or all are vertical.  \\[2mm]
(b) After removing $o(n)$ points from $\pts$, we also have: \begin{itemize}[noitemsep,topsep=1pt,leftmargin=5mm]
    \item Every line $\ell\in \lines$ satisfies $|\ell\cap \pts| = \Theta(\sqrt{n})$. When thinking of $\ell$ as $\RR$, the points of $\ell\cap \pts$ are contained in a generalized arithmetic progression of dimension $\Theta(1)$ and size $\Theta(\sqrt{n})$.
    \item There exist generalized arithmetic progressions $A_1,\ldots,A_s$ with $s=\Theta(1)$, each of constant dimension and size $O(\sqrt{n})$. For every $\ell\in\lines$, there exists $r\in \RR$ and $1\le j \le s$, such that $|(\ell\cap \pts) \cap (r+A_j)| = \Theta(\sqrt{n})$. (When thinking of $\ell$ as $\RR$.)
    \item The lines of $\lines$ can be partitioned into $\Theta(1)$ disjoint subsets, each of size $\Theta(\sqrt{n})$. In each subset, the intercepts of the lines (with the orthogonal axis) are contained in a generalized arithmetic progression of dimension $\Theta(1)$ and size $\Theta(\sqrt{n})$.  
    \end{itemize}
\end{theorem}

When seeing Theorem \ref{th:Structure}, one might conjecture that the point set must be contained in a Cartesian product of a relatively small size. 
In other words, that a set with $\Theta(\sqrt{n})$ distinct distances can be covered by $\Theta(\sqrt{n})$ vertical lines and also by $\Theta(\sqrt{n})$ horizontal lines.  
Surprisingly, this is far from the case. Sets with $\Theta(\sqrt{n})$ distinct distances may be very different than Cartesian products.

\begin{theorem} \label{th:Construction}
There exists a set $\pts$ of $n$ points such that no two points of $\pts$ have the same $x$-coordinate and $D(\pts)=2\sqrt{n}-2$ under the $\ell_\infty$ metric.
\end{theorem}

Instead of the structure obtained in Theorem \ref{th:SzekelyImprov}, one may find a smaller subset with more structure. 

\begin{corollary}\label{co:structureSubset}
Let $\pts$ be a set of $n$ points such that $D(\pts)=\Theta(\sqrt{n})$ under the $\ell_\infty$ metric.
Then there exist a subset $\pts' \subseteq \pts$ of $\Theta(n)$ points with the following properties:
\begin{itemize}[noitemsep,topsep=1pt,leftmargin=5mm]
    \item There exists a set $\lines$ of $\Theta(\sqrt{n})$ lines, either all vertical or all horizontal. Every line $\ell\in \lines$ satisfies $|\ell\cap \pts'| = \Theta(\sqrt{n})$.
    \item There exists a generalized arithmetic progression $A$ with of dimension $\Theta(1)$ and size $\Theta(\sqrt{n})$. For every $\ell\in\lines$, when thinking of $\ell$ as $\RR$, there exists $r\in \RR$ such that $(\ell\cap \pts') \subseteq (r+A)$.
    \item The intercepts of the lines (with the orthogonal axis) are contained in another generalized arithmetic progression of dimension $\Theta(1)$ and size $\Theta(\sqrt{n})$.  
    \end{itemize}
\end{corollary}

Theorem \ref{th:SzekelyImprov} is proved in Section \ref{sec:SzekelyImp}.
Theorem \ref{th:Structure}, Theorem \ref{th:Construction}, and Corollary \ref{co:structureSubset} are proved in Section \ref{sec:Structure}.

\section{Strictly convex $\ell_p$ metrics} \label{sec:SzekelyImp}

The goal of this section is to prove Theorem \ref{th:SzekelyImprov}. Throughout the section, we assume that we are working with an $\ell_p$ metric,
for a finite integer $p\ge 3$. 
As stated in the introduction, this result is already known when $p=2$. 
We first present a variety of tools required for the proof of Theorem \ref{th:SzekelyImprov}.

\parag{Bisectors of $\ell_p$ metrics.}
Let $\bis(u,v)$ denote the \emph{bisector} of the distinct points $u,v \in \RR^2$. 
That is, $\bis(u,v)$ is the set of points of $\RR^2$ that are equidistant from $u$ and $v$.
Writing $u=(u_x,u_y)$ and $v=(v_x,v_y)$, the bisector $\bis(u,v)$ is defined by
\begin{equation} \label{eq:BisectDef}
|x-u_x|^p + |y-u_y|^p - |x-v_x|^p - |y-v_y|^p = 0. 
\end{equation}

\begin{figure}[h]
    \centering
    \includegraphics[scale=0.35]{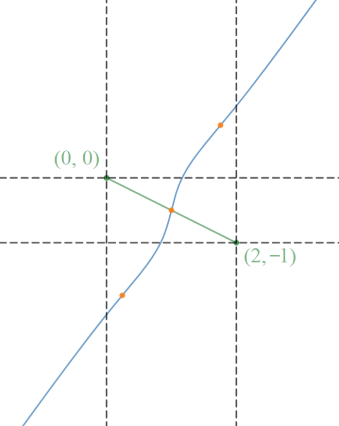}
    \caption{\small \sf In each of the nine regions, the bisector is a segment of an algebraic curve of degree at most $p$. Each inflection point is in a different region that contains a bounded segment of the bisector.} 
    \label{fi:L_pBis-Inflx}
\end{figure}

When $p$ is even, we may remove the absolute values from \eqref{eq:BisectDef}.
We then get that $\bis(u,v)$ is an algebraic curve of degree at most $p-1$.\footnote{An algebraic curve $\gamma\subset \RR^2$ is a one-dimensional set that is the zero set of a polynomial of $\RR[x,y]$. The degree of $\gamma$ is the minimum degree over all polynomials whose zero set is $\gamma$.}
When $p$ is odd, we partition $\RR^2$ into at most nine regions with the lines $x=u_x$, $x=v_x$, $y=u_y$, and $y=v_y$. In each region, we may replace the absolute values in \eqref{eq:BisectDef} with specific signs.
Thus, the intersection of $\bis(u,v)$ with each region is a segment of an algebraic curve of degree at most $p$.
See Figure \ref{fi:L_pBis-Inflx}.

Garibaldi \cite{Garibaldi04} studied the behavior of bisectors under strictly convex $\ell_p$ metrics. 
We now describe some of her results. 

\begin{lemma} \label{le:BisectorProps}
Consider two distinct points $u,v\in \RR^2$. \\[2mm]
(a) The bisector $\bis(u,v)$ is a line when the line containing $u$ and $v$ has slope $0,1,-1,$ or $\infty$. Otherwise, $\bis(u,v)$ contains no line segment.\\
(b) The bisector $\bis(u,v)$ is monotone. \\
(c) Cut $\bis(u,v)$ into two pieces at the midpoint of $u$ and $v$. These two pieces are identical up to a rotation. \\
(d) If $\bis(u,v)$ is not a line and $\bis(u,v) = \bis(u',v')$, then $(u,v)=(u',v')$. \\
(e) There exists an absolute constant $C$ such that every two non-identical bisectors intersect in at most $C$ points. 
\end{lemma}

Since a non-line bisector $\bis(u,v)$ is monotone and does not contain line segments, it is a graph function with respect to either axis. That is, we can define $\bis(u,v)$ both as $y=f(x)$ and as $x=g(y)$.  
We say that a point $(u_x,u_y)$ is an \emph{inflection point} of $\bis(u,v)$ if $f''(x)=0$ or $g''(y)=0$.
The following results are part of the proof of Lemma 5.4.1 in Garibaldi \cite{Garibaldi04}.

\begin{lemma}\label{le:BisectorProps2}
Consider two distinct points $u=(u_x,u_y)$ and $v=(v_x,v_y)$, such that $\bis(u,v)$ is not a line. Partition $\RR^2$ into regions using the lines $x=u_x$, $x=v_x$, $y=u_y$, and $y=v_y$. \\[2mm]
(a) The bisector $\bis(u,v)$ has exactly three inflection points. One of these is the midpoint of $u$ and $v$.\\
(b)  The bisector $\bis(u,v)$ intersects five of the regions.
Two of these regions contain unbounded segments of $\bis(u,v)$. Each of the other three regions contains an inflection point. (See Figure \ref{fi:L_pBis-Inflx}.) 
\end{lemma}

\parag{Crossings in multigraphs.}
Two edges in a graph are said to be \emph{parallel} if they have the same endpoints.
A \emph{mutligraph} is a graph that may contain parallel edges.
The \emph{multiplicity} of an edge $e$ is the number of edges with the same endpoints as $e$ (including $e$ itself).

In a \emph{drawing} of a graph, every vertex is a distinct point in the plane and
every edge is a Jordan arc connecting the two corresponding vertices. We assume that the interior
of every such arc does not contain vertices, that any two arcs have a finite number of intersections,
and that no three arcs intersect at the same point. (This can be achieved while maintaining the crossings, by slightly perturbing the drawing.) The crossing number $\crs(G)$ of a graph
$G$ is the minimum number of edge crossings across all drawings of $G$.

\begin{lemma}[Crossing lemma for multigraphs \cite{Szek97}] \label{le:crossing} Consider a multigraph $G = (V, E)$ with $|V| = n$, $|E| = e$, and maximum edge multiplicity $m$. If $e > 5mn$, then
\[ \crs(G) = \Omega\left(\frac{e^3}{mn^2}\right).\] \end{lemma}

\parag{Incidences with curves.}
One can define an infinite family of algebraic curves in $\RR^2$ with a polynomial $f\in \RR[x,y]$ whose coefficients depend on parameters $s_1,\ldots,s_k\in \RR$.
Each curve in the family is obtained by assigning values to the parameters and then taking the zero set of $f$.
For example, the family of all circles can be defined with the polynomial $(x-c_x)^2+(y-c_y)^2 -r^2$ and parameters $c_x,c_y,r\in \RR$.
The \emph{dimension} of such a family of curves is the number of parameters.  

Let $\pts$ be a set of points and let $\curves$ be a set of curves, both in $\RR^2$.
An \emph{incidence} is a pair $(p,\gamma)\in\pts\times \curves$ such that the point $p$ is on the curve $\gamma$.
We denote the number of incidences in $\pts\times\curves$ as $I(\pts,\curves)$.
The following incidence result is by Sharir and Zahl \cite{SZ16}. 

\begin{theorem} \label{th:SharirZahl}
Let $\pts$ be a set of $m$ points.
Let $\curves$ be a set of $n$ curves that belong to an $s$-dimensional family of curves of degree at most $D$, no two sharing an irreducible component. Then, for every $\eps>0$,
$$I(\pts,\curves)=O_{s,D,\eps}\left(m^{\frac{2s}{5s-4}}n^{\frac{5s-6}{5s-4}+\eps} + m^{2/3}n^{2/3}+m+n\right).$$
\end{theorem}

Consider a point set $\pts\subset \RR^2$ and $r\in\RR$. A curve $\gamma$ is \emph{$r$-rich with respect to} $\pts$ if $|\gamma\cap\pts|\ge r$.
When the point set is clear from the context, we simply write that $\gamma$ is $r$\emph{-rich}. 

We are now ready to prove Theorem \ref{th:SzekelyImprov}. 
We first restate this theorem.
\vspace{2mm}

\noindent {\bf Theorem \ref{th:SzekelyImprov}}
\emph{Let $\pts$ be a set of $n$ points in $\RR^2$ and let $p> 1$ be an integer. Then, under the $\ell_p$ metric and for any $\eps>0$, there exists a point $u\in \pts$ such that} 
\[ D(u,\pts) = \Omega_{p,\eps}\left(n^{6/7-\eps}\right).\]

\begin{proof}
We adapt Sz\'ekely's proof \cite{Szek97} for Euclidean metric. 
We may assume that $p\ge 3$, since the case of $p=2$ was handled in \cite{SoTo01}.
Set $t=\max_{u\in \pts} D(u,\pts)$. 
In other words, $t$ is the maximum number of distances between a point of $\pts$ and the rest of $\pts$. We may assume that $t=O(n/\log^2 n)$, since otherwise we are done. 

Throughout the proof, when mentioning circles, we refer to circles defined by the $\ell_p$ norm (rather than to Euclidean circles). That is, the set of points at a given distance from a given point. 
Let $\circs$ be the set of circles centered at a point of $\pts$ and incident to another point of $\pts$. By definition, each point of $\pts$ has at most $t$ such circles centered at it. This implies that $|\circs|\le nt$. 
Since each point of $\pts$ is incident to a circle centered at every other point, we have that $I(\pts,\circs) = n(n-1)$. We remove from $\circs$ all circles incident to at most two points of $\pts$. 
This decreases the number of incidences by at most $2nt$, so we still have that $I(\pts,\circs) = \Theta(n^2)$. 

Consider the graph $G=(V,E)$, defined as follows. The set $V$ contains a vertex for each point of $\pts$. 
For every circle $\gamma\in\circs$, we add to $E$ an edge between every two points that are consecutive along $\gamma$. 
Note that pairs of points that are consecutive in multiple circles of $\circs$ lead to parallel edges in $G$. 
The proof proceeds by double counting $\crs(G)$.

We draw every vertex of $V$ as its corresponding point and every edge of $E$ as the corresponding circular arc from $\circs$.
In this drawing, edges intersect only at the intersection points of the circles of $\circs$. Two circles have at most two intersection points (for example, see \cite[Lemma 1.4]{IL04}). This implies 
\begin{equation} \label{eq:CrossingUpper} 
\crs(G) \le 2\cdot \binom{|\circs|}{2} \leq 2\cdot \binom{nt}{2} =\Theta(n^2t^2). 
\end{equation}

A circle of $\circs$ that is incident to $j\ge 3$ points of $\pts$ contributes $j$ edges to $E$. Since $I(\pts,\circs)=\Theta(n^2)$, we get that $|E|=\Theta(n^2)$. To get a lower bound for $\crs(G)$, we wish to rely on Lemma \ref{le:crossing}. However, since some edges in $E$ may have a high multiplicity, this leads to a weak bound. Thus, we first study edges with high multiplicity. 

Consider distinct points $u,v\in \pts$.
Every edge in $E$ between $u$ and $v$ corresponds to a distinct point of $\pts$ that is equidistant from $u$ and $v$.
In other words, every edge between $u$ and $v$ corresponds to a point of $\pts$ that is incident to $\bis(u,v)$.
Thus, edges with high multiplicity correspond to bisectors that contain many points of $\pts$.
For $k$ whose value will be determined below, we study edges with multiplicity at least $k$.
We separately handle bisectors that are lines and other bisectors.

\parag{Line bisectors.} 
By Lemma \ref{le:BisectorProps}, every line bisector has slope $0,1,-1,$ or $\infty$.
Fix an integer $j\ge 0$. Since every point is incident to at most one line of slope 1, the number $2^j$-rich lines of slopes 1 is at most $n/2^j$. 
Similarly, the number of $2^j$-rich lines with one of the four possible slopes is $O(n/2^j)$. 

Let $T_j\subset \pts^3$ be the set of triples $(u,v,w)$ where the multiplicity of the edge $(u,v)$ is at least $k$, the bisector $\bis(u,v)$ is $2^j$-rich, and $w\in \bis(u,v)$. Note that every such triple corresponds to at least $2^j$ edges with multiplicity at least $k$.
We now derive an upper bound for $|T_j|$.

Let $\ell$ be a $2^j$-rich line with one of the four aforementioned slopes. 
Let $w\in \pts$ be a point incident to $\ell$.
Since the circles of the $\ell_p$ metric are strictly convex, every circle centered at $w$ intersects $\ell$ twice.
By definition, at most $t$ such circles are incident to points of $\pts$. Each circle contains at most two pairs $(u,v)\in \pts^2$ that are consecutive along the circle and satisfy $\ell=\bis(u,v)$. 
We conclude that $|T_j| = O(tn/k)$.

We partition the edges of $E$ according to their multiplicity. 
Specifically, for each $\log k \le j \le \log n$, we consider edges with multiplicity at least $2^j$ and smaller than $2^{j+1}$. This implies that the number of edges with a line bisector and multiplicity at least $k$ is smaller than
\[ \sum_{j=\log k}^{\log n} |T_j| \cdot 2^{j+1} =\sum_{j=\log k}^{\log n} O\left(\frac{tn}{2^{j}} \cdot 2^{j+1}\right)= tn\sum_{j=\log k}^{\log n}O(1) = O(tn\log n). \]

We remove the above edges.
The assumption $t=O(n/\log^2 n)$ implies that the number of removed edges is $O(n^2/\log n)$. We thus still have that $|E|=\Theta(n^2)$.

\parag{Non-line bisectors.}
We first consider the case where $p$ is even.
In this case, the non-line bisectors are algebraic curves of degree at most $p-1$.
Fix an integer $j$ that is at least some sufficiently large constant. Let $\curves$ be the set of $2^j$-rich non-line bisectors of pairs of points from $\pts$. By Lemma \ref{le:BisectorProps}, the bisectors of $\curves$ are distinct and no two bisectors  share an irreducible component.  
Since each bisector is defined by two planar points, $\curves$ is part of a four-dimensional family of algebraic curves. 
For any $\eps'>0$, Theorem \ref{th:SharirZahl} implies 
\begin{equation} \label{eq:SharirZahlApplication}
I(\pts,\curves)=O_{p,\eps'}\left(n^{1/2}|\curves|^{7/8+\eps'}+n^{2/3}|\curves|^{2/3}+n+|\curves|\right).
\end{equation}

Since every curve of $\curves$ is incident to at least $2^j$ points of $\pts$, we have that  $I(\pts,\curves)\geq 2^j|\curves|$.
Combining this with the above upper bound and setting $\eps'' = 32\eps'/(1-8\eps')$ yields 
\begin{equation} \label{eq:richBisect} |\curves| = O_{p,\eps'}\left(\frac{n^{4/(1-8\eps')}}{2^{8j/(1-8\eps')}}+\frac{n^2}{2^{3j}}+\frac{n}{2^j}\right) = O_{p,\eps'}\left(\frac{n^{4+\eps''}}{2^{8j}}+\frac{n^2}{2^{3j}}+\frac{n}{2^j}\right). \end{equation}
(Since $j$ is assumed to be sufficiently large, the fourth term in \eqref{eq:SharirZahlApplication} cannot dominate the bound.)

The first term on the right-hand side of \eqref{eq:richBisect} dominates the bound when $2^j=O(n^{(2+ \eps'')/5})$. 
The second term of \eqref{eq:richBisect} dominates when $2^j=\Omega(n^{(2+ \eps'')/5})$ and $2^j =O(\sqrt{n})$. 
The last term of \eqref{eq:richBisect} dominates when $2^j =\Omega(\sqrt{n})$. 
As in the case of line bisectors, we dyadically decompose the edges of $E$ according to their multiplicity. 
In this case, each bisector corresponds to at most one pair $(u,v)\in \pts^2$.
Thus, the number of edges with a non-line bisector and multiplicity at least $k$ is 
\begin{align*}
O_{p,\eps'}&\left(\sum_{j=\log k}^{\frac{2+\eps''}{5} \log n} \frac{n^{4+\eps''}}{2^{8j}}\cdot 2^{j+1} + \sum_{j=\frac{2+\eps''}{5} \log n}^{\frac{1}{2} \log n} \frac{n^2}{2^{3j}}\cdot 2^{j+1} +
\sum_{j=\frac{1}{2} \log n}^{\log n} \frac{n}{2^{j}}\cdot 2^{j+1}\right)
\\[2mm] & \hspace{20mm}= O_{p,\eps'}\left(n^{4+\eps''}\cdot \sum_{j=\log k}^{\frac{2+\eps''}{5} \log n} \frac{1}{2^{7i}} + n^2\cdot \sum_{j=\frac{2+\eps''}{5} \log n} \frac{1}{2^{2j}}+
n \cdot\sum_{j=\frac{1}{2} \log n}^{\log n} 1\right)
\\[2mm] & \hspace{20mm}= O_{p,\eps'}\left(\frac{n^{4+\eps''}}{k^7} +
n^{\frac{6-2\eps''}{5}} + n\log n\right).
\end{align*}

Taking $k=C_{p,\eps'}n^{(2+\eps'')/7}$, for an appropriate constant $C_{p,\eps'}$, leads to at most $n^2/100$ edges with multiplicity at least $k$ and a non-line bisector.
As in the line bisector case, We remove those edges, while still having that $|E|=\Theta(n^2)$.

\parag{Completing the even case.}
After the above pruning of $E$, every remaining edge has multiplicity smaller than $C_{p,\eps'}n^{(2+\eps'')/7}$.
We may thus apply Lemma \ref{le:crossing} with $m=C_{p,\eps'}n^{(2+\eps'')/7}$, obtaining
\[ \crs(G) = \Omega\left(\frac{n^6}{n^{(2+\eps'')/7}\cdot n^2}\right) = \Omega\left(n^{\frac{26-\eps''}{7}}\right). \]

Combining this with \eqref{eq:CrossingUpper} leads to $n^{\frac{26-\eps''}{7}}=O(n^2t^2)$.
Simplifying gives $t=\Omega\left(n^{\frac{6}{7}-\frac{\eps''}{14}}\right)$.
To complete the theorem, we choose $\eps'$ such that $\eps = \eps''/14$.

\parag{The case of odd $p$.}
We wish to repeat the proof of the even case.
The line bisectors are handled as in the even case.
However, the non-line bisectors are not algebraic curves, since their defining equations include absolute values.
By Lemma \ref{le:BisectorProps2}, every non-line bisector can be cut into five segments, such that each segment is contained in an algebraic curve of degree at most $p$.

We wish to apply Theorem \ref{th:SharirZahl}, but we have segments of algebraic curves rather than algebraic curves.
We thus replace each segment with the algebraic curve that contains it.\footnote{More formally, we take the \emph{Zariski closure} of each segment.}
This leads to a new potential issue: Segments originating from different bisectors may be contained in the same irreducible algebraic curve.

From \eqref{eq:BisectDef}, the curves that contain the segments are all of the form 
\begin{equation} \label{eq:containingCurve}
(x-a_1)^p \pm (x-a_2)^p = \pm (y-a_3)^p \pm (y-a_4)^p, 
\end{equation}
where $a_1,a_2,a_3,a_4\in \RR$. Each $\pm$ symbol may represent a different sign.  

By Lemma \ref{le:BisectorProps2}, every bisector segment is either unbounded or contains an inflection point.
By inspecting \eqref{eq:containingCurve}, we note that such a curve may contain at most four unbounded segments. Indeed, unbounded segments may only occur when $x,y\to \pm \infty$. 

Recall that inflection points occur when $\frac{\partial^2}{\partial^2 x} y =0$ and when $\frac{\partial^2}{\partial^2 y} x =0$. We now study the former, and the latter can be handled symmetrically. We take a derivative of \eqref{eq:containingCurve} according to $x$ and rearrange, obtaining
\[ \frac{\partial}{\partial x} y = \frac{(x-a_1)^{p-1}\pm(x-a_2)^{p-1}}{\pm(y-a_3)^{p-1}\pm(y-a_4)^{p-1}}.\]
Taking another derivative leads to
\[ \frac{\partial^2}{\partial^2 x} y = \frac{f(x,y)}{(\pm(y-a_3)^{p-1}\pm(y-a_4)^{p-1})^3},\]
where $f\in \RR[x,y]$ is of degree smaller than $3p$.

The inflection points of the curve are the solutions to \eqref{eq:containingCurve} that also satisfy $f(x,y)=0$.
Note that \eqref{eq:containingCurve} and the set of solutions to $f(x,y)=0$ cannot share irreducible components, since no non-line irreducible curve consists entirely of inflection points.
By B\'ezout's theorem (for example, see \cite{CLOu}), there are fewer than $3p^2$ solutions to this system.
Including also four unbounded segments and derivatives according to $y$, the number of bisector segments contained in a curve of the form \eqref{eq:containingCurve} is smaller than $6p^2+4$.

We apply Theorem \ref{th:SharirZahl} while keeping one copy of each curve that repeats multiple times. To obtain a valid upper bound for the number of incidences, we multiply the bound of the theorem by  $6p^2+4$. 
The rest of the proof is identical to that of the even case.
\end{proof}

\section{Structure in $\ell_1$ and $\ell_\infty$} \label{sec:Structure}

Our goal in this section is to prove Theorem \ref{th:Structure}, Theorem \ref{th:Construction}, and Corollary \ref{co:structureSubset}. 
While these results are stated for the $\ell_\infty$ metric, they equally apply to the $\ell_1$ metric. 
Indeed, the unit circle of the $\ell_\infty$ metric is an axis-parallel square of side length 2 and the unit circle of the $\ell_1$ metric is a square rotated by $\pi/4$ of side length $\sqrt{2}$. 
For a point set $\pts$, let $\pts'$ be the set obtained from $\pts$ after a rotation of $\RR^2$ by an angle of $\pi/4$ followed by a uniform scaling by a factor of $\sqrt{2}$.
Then the set of distances spanned by $\pts$ with the $\ell_1$ metric is identical to the set of distances spanned by $\pts'$ with the $\ell_\infty$ metric. A symmetric argument holds when moving from $\ell_\infty$ to $\ell_1$.

Throughout this section, we work with the $\ell_\infty$ metric.
To prove our results, we first need to introduce some basic tools and concepts from additive combinatorics. For reference and additional details, see Tao and Vu \cite{TV06}.
The \emph{difference set} of a set $A\subset \RR$ is
\[ A-A = \{a-a'\ :\ a,a'\in A\}. \]

When $|A|=n$, we have that $2n-1\le |A-A| \le \binom{n+1}{2}$.
The only sets of $n$ numbers that satisfy $|A-A|=2n-1$ are arithmetic progressions. 
However, sets that satisfy $|A-A|=\Theta(n)$ have a significantly more varied behavior. 
A similar situation happens for sets of $n$ points in $\RR^2$ that span few distinct distances: The only sets that span $\sqrt{n}-1$ distinct distances are $\sqrt{n}\times\sqrt{n}$ uniform lattices. However, sets that span $\Theta(\sqrt{n})$ distinct distances can be rather different.

A \emph{generalized arithmetic progression} of \emph{dimension} $d$ is defined as
\[ \left\{ a + \sum_{j=1}^{d} k_j b_j :\, a,b_1,\ldots,b_d \in \RR \text{ and with integer } 0\le k_j \le n_j-1 \text{ for every } 1\le j \le d \right\}. \]

The following is a special case of Freiman's theorem. 
\begin{theorem}[Freiman's theorem over $\RR$] \label{th:Freiman}
Let $A\subset \RR$ be a finite set with $|A-A| \le k |A|$ for some constant $k$.
Then $A$ is contained in a generalized arithmetic progression of size at most $cn$ and dimension at most $d$. Both $c$ and $d$ depend on $k$ but not on $|A|$.
\end{theorem}

The \emph{difference energy} of a set $A\subset \RR$ is
\[ E(A) = \left|\{(a_1,a_2,a_3,a_4)\in A^4\ :\ a_1-a_2 = a_3-a_4\}\right|. \]
Difference energy is a main tool in studying additive properties of sets. 
For our purposes, we only need one property of this energy: the Balog--Szemer\'edi--Gowers theorem. The following variant of this theorem is by Schoen \cite{Schoen15}.

\begin{theorem} \label{th:schoen}
Let $A\subset\RR$ such that $E(A)=\delta |A|^3$. Then, there exists $A'\subset A$ such that $|A'|= \Omega(\delta|A|)$ and $|A'-A'|=O(\delta^{-4}|A'|)$.
\end{theorem}

For a set $A\subset \RR$ and $r\in \RR$, we define the \emph{translation} 
\[ r+A = \{a+r\ :\ a:\in A\}.\]

We are now ready to prove Theorem \ref{th:Structure}. We first restate this result.
\vspace{2mm}

\noindent {\bf Theorem \ref{th:Structure}} 
\emph{Let $\pts$ be a set of $n$ points such that $D(\pts)=\Theta(\sqrt{n})$ under the $\ell_\infty$ metric. Then: \\[2mm]
(a) There exists a set $\lines$ of $\Theta(\sqrt{n})$ lines such that $\pts \subset \bigcup_{\ell\in \lines} \ell$. Either all lines of $\lines$ are horizontal, or all are vertical.  \\[2mm]
(b) After removing $o(n)$ points from $\pts$, we also have: \begin{itemize}[noitemsep,topsep=1pt,leftmargin=5mm]
    \item Every line $\ell\in \lines$ satisfies $|\ell\cap \pts| = \Theta(\sqrt{n})$. When thinking of $\ell$ as $\RR$, the points of $\ell\cap \pts$ are contained in a generalized arithmetic progression of dimension $\Theta(1)$ and size $\Theta(\sqrt{n})$.
    \item There exist generalized arithmetic progressions $A_1,\ldots,A_s$ with $s=\Theta(1)$, each of constant dimension and size $O(\sqrt{n})$. For every $\ell\in\lines$, there exists $r\in \RR$ and $1\le j \le s$, such that $|(\ell\cap \pts) \cap (r+A_j)| = \Theta(\sqrt{n})$. (When thinking of $\ell$ as $\RR$.)
    \item The lines of $\lines$ can be partitioned into $\Theta(1)$ disjoint subsets, each of size $\Theta(\sqrt{n})$. In each subset, the intercepts of the lines (with the orthogonal axis) are contained in a generalized arithmetic progression of dimension $\Theta(1)$ and size $\Theta(\sqrt{n})$.  
    \end{itemize}
}    
\begin{proof}
(a) Let $(x_1, y_1)$ be the point of $\pts$ that minimizes the value of $x+y$. Let $(x_2, y_2)$ be the point that maximizes this value. Set $I_1^+ = x_1+y_1$ and $I_4^+=x_2+y_2$. The meaning of this notation will become clear below. For now, note that each $I^+$ can be thought of as the projection of a point on the line $y=-x$. There might be multiple points satisfying $x+y=I_1^+$ and $x+y=I_4^+$. Among those, we choose $(x_1, y_1)$ and $(x_2, y_2)$ that minimize $(y_1-x_1)-(y_2-x_2)$.
Let $(x_3, y_3)$ be the point of $\pts$ that minimizes the value of $y-x$. Let $(x_4, y_4)$ be the point that maximizes this value. Set $I_1^- = y_3-x_3$ and $I_4^-=y_4-x_4$.
When there are multiple points that lead to these values, We choose $(x_3, y_3)$ and $(x_4, y_4)$ that minimize $(x_3+y_3)-(x_4+y_4)$.

By definition, $\pts$ is contained in the rectangle $R$ defined by $I_1^+\le x+y \le I_4^+$ and $I_1^- \le y-x\le I_4^-$. We cut $R$ into at most nine rectangular pieces using:  
\begin{itemize}[noitemsep,topsep=1pt] 
    \item The line of slope 1 incident to $(x_1, y_1)$.
    \item The line of slope 1 incident to $(x_2, y_2)$. 
    \item The line of slope -1 incident to $(x_3, y_3)$.
    \item The line of slope -1 incident to $(x_4, y_4)$.
\end{itemize}
An example is depicted in Figure \ref{fi:NineRects}.
There are fewer than nine rectangles  when two of the lines are identical or when a line contains a side of $R$. 

Set $\pts' = \{(x_1, y_1),(x_2, y_2),(x_3, y_3),(x_4, y_4)\}$. 
These four points are not necessarily distinct.
If a point appears more than once in $\pts'$, then we only keep one copy of that point.
That is, $\pts'$ is a set of at most four distinct points. 

We consider the case where $y_1-x_1 \ge y_2-x_2$. The case where $y_1-x_1 < y_2-x_2$ can be handled symmetrically.
We set $I_3^-=y_1-x_1$ and $I_2^- = y_2-x_2$.
Note that $I_1^- \le I_2^- \le I_3^- \le I_4^-$.
The analysis of the case of $y_1-x_1 \ge y_2-x_2$ is further divided into two cases, as follows.

\parag{Case 1.} In this case, we assume that $x_3+y_3\le x_4+y_4$. 
We set $I_3^+=x_4+y_4$ and $I_2^+ = x_3+y_3$.
Note that $I_1^+ \le I_2^+ \le I_3^+ \le I_4^+$.
We define nine (possibly not distinct) rectangles:
\begin{align*}
R_1 &= \{(a,b)\ :\ I_1^- \leq b-a \leq I_2^-,\ I_1^+ \leq a+b \leq I_2^+\}, \\
R_2 &= \{(a,b)\ :\ I_1^- \leq b-a \leq I_2^-,\ I_2^+ < a+b < I_3^+ \}, \\
R_3 &= \{(a,b)\ :\ I_1^- \leq b-a \leq I_2^-,\ I_3^+ \le a+b \le I_4^+ \}, \\
R_4 &= \{(a,b)\ :\ I_2^- < b-a < I_3^-,\ I_1^+ \le a+b \le I_2^+ \}, \\
R_4 &= \{(a,b)\ :\ I_2^- < b-a < I_3^-,\ I_1^+ \le a+b \le I_2^+ \}, \\
R_5 &= \{(a,b)\ :\ I_2^- < b-a < I_3^-,\ I_2^+ < a+b < I_3^+ \}, \\
R_6 &= \{(a,b)\ :\ I_3^- \le b-a \le I_4^-,\ I_2^+ < a+b < I_3^+ \}, \\
R_7 &= \{(a,b)\ :\ I_1^- \le b-a \le I_2^-,\ I_3^+ \le a+b \le I_4^+ \}, \\
R_8 &= \{(a,b)\ :\ I_2^- < b-a < I_3^-,\ I_3^+ \le a+b \le I_4^+ \}, \\
R_9 &= \{(a,b)\ :\ I_3^- \le b-a \le I_4^-,\ I_3^+ \le a+b \le I_4^+ \}.
\end{align*}

Figure \ref{fi:NineRects} depicts Case 1 when there are nine distinct rectangles.

\begin{figure}[htp]
    \centering
    \includegraphics[scale=0.4]{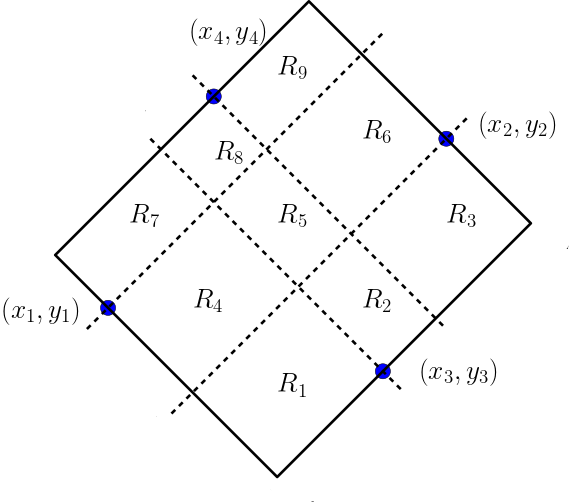}
    \caption{\small \sf Case 1 when there are nine distinct rectangles.}
    \label{fi:NineRects}
\end{figure}
 
Consider a point $u \in \pts'$. Since $D(\pts)=\Theta(\sqrt{n})$, the points of $\pts\setminus \{u\}$ are contained in a set ${\mathcal Q}_u$ of $O(\sqrt{n})$ axis-parallel squares centered at $u$. 
By the pigeonhole principle, there exists a square in ${\mathcal Q}_u$ that contains at least $n/|{\mathcal Q}_p|$ points of $\pts$. It is not difficult to check that these points span $\Omega(n/|{\mathcal Q}_u|)$ distinct distances. Since
$D(\pts)=\Theta(\sqrt{n})$, we get that $|{\mathcal Q}_u|= \Omega(\sqrt{n})$. 
That is, $|{\mathcal Q}_u|= \Theta(\sqrt{n})$.

Let ${\mathcal Q}$ be the set of $\Theta(\sqrt{n})$ squares that are centered at a point of $\pts'$ and contain at least one point of $\pts$ on their boundary.
Let $\lines_{H}$ be the set of $\Theta(\sqrt{n})$ horizontal lines that contain the side of at least one square of ${\mathcal Q}$.
Let $\lines_{V}$ be the set of $\Theta(\sqrt{n})$ vertical lines that contain the side of at least one square of ${\mathcal Q}$.
Note that every point of $\pts$ is contained in at least one line of $\lines_{H}\cup \lines_{V}$.
 
Consider a point $u\in \pts\setminus \pts'$ that lies in the rectangle $R_1$.
Let $q_1 \in {\mathcal Q}$ be the square centered at $(x_1,y_1)$ and containing $u$ on its boundary. 
Note that $u$ must be on a vertical side of $q_1$.
Similarly, when considering a square $q_3 \in {\mathcal Q}$ centered at $(x_3,y_3)$, the point $u$ is also on a vertical side of $q_3$. 
When considering squares  $q_2,q_4 \in {\mathcal Q}$ centered at $(x_2,y_2)$ and $(x_4,y_4)$ respectively, the point $u$ is on a horizontal side.
   
An analysis similar to the one in the preceding paragraph  holds for $u \in \pts\setminus \pts'$ in any rectangle $R_{i}$. We associate with each rectangle a quadruple $(\alpha_1,\alpha_2,\alpha_3,\alpha_4)\in \{V,H\}^4$, where $\alpha_j=V$ means that a vertical side of a square around $(x_j,y_j)$ has $u$ on its boundary, and $\alpha_j=H$ means a horizontal side. 
For example, in the preceding paragraph, $R_1$ is associated with the quadruple $(V,H,V,H)$.
More generally: 
\begin{align*}
&R_1:(V,H,V,H), \hspace{11mm}  R_2: (V,H,H,H), \hspace{11mm} R_3: (V, H, H, V),  \\
&R_4: (V, V, V, H), \qquad \quad \: R_5: (V, V, H, H), \qquad \quad \: R_6: (V, V, H, V),  \\
&R_7: (H, V, V, H), \hspace{11mm} R_8: (H, V, H, H), \hspace{11mm} R_9: (H, V, H, V). 
\end{align*}

Consider a point $u\in \pts\setminus \pts'$. Note that, no matter what rectangle $u$ is in, it is on at least one horizontal side of a square and at least one vertical side of a square. In other words, $u$ is incident to a line of $\lines_{H}$ and to a line of $\lines_{V}$. Part (a) of the theorem is obtained by setting $\lines$ to be either $\lines_{H}$ or $\lines_{V}$. That is, in this case we may choose the direction of the lines of $\lines$.
 
\parag{Case 2.} In this case, we assume that $x_3+y_3> x_4+y_4$. 
We set $I_3^+=x_3+y_3$ and $I_2^+ = x_4+y_4$.
As before, $I_1^+ \le I_2^+ \le I_3^+ \le I_4^+$.
We define the nine rectangles $R_1,\ldots,R_9$ as in Case 1, but with the new values of $I_2^+$ and $I_3^+$.
Figure \ref{fi:NineRects2} depicts Case 2 when there are nine distinct rectangles.

\begin{figure}[htp]
    \centering
    \includegraphics[scale=0.43]{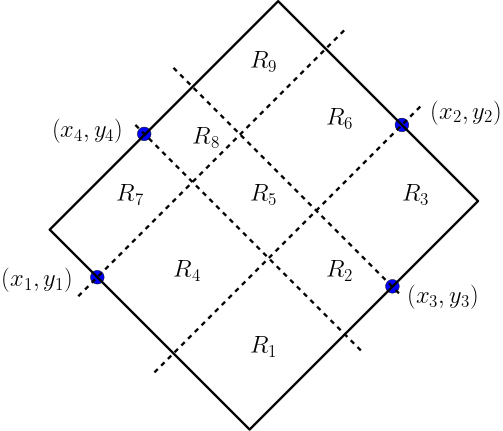}
    \caption{\small \sf Case 2 when there are nine distinct rectangles.}
    \label{fi:NineRects2}
\end{figure}

We define the quadruples $(\alpha_1,\alpha_2,\alpha_3,\alpha_4)\in \{V,H\}^4$ as in Case 1.
In this case, we have the quadruples
\begin{align*}
&R_1:(V, H, V, H), \hspace{11mm}  R_2: (V, H, V, V), \hspace{11mm} R_3: (V, H, H, V),  \\
&R_4: (V, V, V, H), \qquad  \quad R_5: (V, V, V, V), \qquad \quad \: R_6: (V, V, H, V),  \\
&R_7: (H, V, V, H), \hspace{11mm} R_8: (H, V, V, V), \hspace{11mm} R_9: (H, V, H, V). 
\end{align*}
 
Unlike Case 1, only vertical lines pass through $R_5$. In the symmetric case when $y_1-x_1 < y_2-x_2$, this cell is crossed only by horizontal lines. Thus, we do not have a choice between $\lines_{H}$ or $\lines_{V}$. We set $\lines$ to be either $\lines_H$ or $\lines_V$, according to the slopes in $R_5$.  
\vspace{2mm}

(b) Since $D(\pts)=\Theta(\sqrt{n})$, every line of $\lines$ contains  $O(\sqrt{n})$ points of $\pts$.
We say that a line of $\lines$ is \emph{rich} if it contains $\Theta(\sqrt{n})$ points of $\pts$. 
Since there are $O(\sqrt{n})$ non-rich lines in $\lines$, these lines contain $o(n)$ points of $\pts$.
In other words, we may assume that the number of points of $\pts$ on non-rich lines is smaller than $\eps n$, for any $\eps>0$.
This implies that $(1-\eps)n$ points of $\pts$ are on the rich lines. Thus, there are $\Theta(\sqrt{n})$ rich lines.

We discard the non-rich lines from $\lines$. We also discard from $\pts$ the points that were on non-rich lines.
By the above, we remain with $\Theta(\sqrt{n})$ lines in $\lines$ and may ask that $|\pts|\ge n(1-\eps)$ for any $\eps>0$.

Without loss of generality, assume that the lines of $\lines$ are horizontal.
Consider a line $\ell\in \lines$.
By the above pruning, $\ell$ contains $\Theta(\sqrt{n})$ points of $\pts$.
Let $C_\ell$ denote the set of $x$-coordinates of the points of $\ell\cap \pts$.
Note that $C_\ell-C_\ell$ consists of the distances spanned by points of $\ell\cap \pts$, the same distances with a negative sign, and 0. 
This leads to $|C_\ell-C_\ell|=\Theta(\sqrt{n})=\Theta(|C_\ell|)$. Theorem \ref{th:Freiman} implies that $C_\ell$ is contained in a generalized arithmetic progression of constant dimension and size $\Theta(\sqrt{n})$.
We denote this progression as $A_\ell$.

\parag{Studying the arithmetic progressions.}
Fix a line $\ell\in \lines$.
By the above, $C_\ell$ is contained in a generalized arithmetic progression 
\[ A_\ell =  \left\{ a+\sum_{j=1}^{d} k_j b_j\ :\ a,b_1,\ldots,b_d \in \RR \quad \text{and }\quad 0\le k_j \le n_j-1 \right\}. \]
Note that 
\[ A_\ell - A_\ell =\left\{ \sum_{j=1}^{d} k_j b_j\ :\ b_1,\ldots,b_d \in \RR \quad \text{and }\quad 1-n_j\le k_j \le n_j-1 \right\}. \]

Consider another line $\ell'\in \lines$, such that $|(C_{\ell'}-C_{\ell'})\setminus(C_{\ell}-C_{\ell})| = o(\sqrt{n})$.
Since 
\[ (C_{\ell}-C_{\ell})\bigcap(C_{\ell'}-C_{\ell'}) \subset (A_{\ell}-A_{\ell}), \] 
we have that  $|(A_{\ell}-A_{\ell}) \bigcap (C_{\ell'}-C_{\ell'})| = \Theta(\sqrt{n})$. 
Thus, there exists $r\in \RR$ such that  $|(A_{\ell}+r)\bigcap C_{\ell'}|=\Theta(\sqrt{n})$.

We now describe a process for creating the progressions $A_1,\ldots,A_s$. 
We say that a line $\ell\in \lines$ is a \emph{saved line} if in a previous step we chose $A_{\ell}$ to be one of the progressions $A_1,\ldots,A_s$.
We go over the lines of $\lines$ one by one. For the line $\ell$ that we consider in the first step, we set $A_1 = A_\ell$.
For a line $\ell'\in \lines$ considered in any other step, we first check whether $|(C_{\ell'}-C_{\ell'})\setminus(C_{\ell}-C_{\ell})| = o(\sqrt{n})$ for a saved line $\ell$.
If no such saved line exists, then we add $A_{\ell'}$ as a new progression $A_j$ (so $\ell'$ becomes a saved line).
Otherwise, there exists a saved line $\ell$ such that $|(C_{\ell'}-C_{\ell'})\setminus(C_{\ell}-C_{\ell})| = o(\sqrt{n})$.
In this case, there exists $r\in \RR$ such that $|C_{\ell'} \cap A_\ell| = \Theta(\sqrt{n})$.

In the above process, there are $\Theta(1)$ steps where we set a new $A_j$ and a new saved line.
Otherwise, the number of distinct distances would be asymptotically larger than $\sqrt{n}$.
This implies that $s=\Theta(1)$. 
At the end of the above process, for every $\ell\in\lines$ there exists $r\in \RR$ and $1\le j \le s$ such that $|(\ell\cap \pts) \cap (r+A_j)| = \Theta(\sqrt{n})$.

\parag{Studying the distances between the lines.} Below we describe a process for finding $\Theta(\sqrt{n})$ lines of $\lines$ whose $y$-intercepts form a generalized arithmetic progression of  constant dimension and size $\Theta(\sqrt{n})$.
We apply this process to find a set of lines $\lines_1$. We then discard the lines of $\lines_1$ from $\lines$ and the corresponding points from $\pts$.
If there remain $\Theta(\sqrt{n})$ lines in $\lines$, then we repeat the process to obtain another set $\lines_2 \subset \lines \setminus \lines_1$.
We then discard the lines of $\lines_2$ and the corresponding points of $\pts$.
We repeat this process until the number of lines that remain in $\lines$ is $o(\sqrt{n})$.

We say that a pair of distinct points $(u,v)\in \pts^2$ is \emph{horizontal} if the distance between $u$ and $v$ is the difference of their $x$-coordinates.
We say that $(u,v)\in \pts^2$ is \emph{vertical} if the distance is the difference of the $y$-coordinates. 
Note that a pair might be both horizontal and vertical.
As long as $\Theta(\sqrt{n})$ lines remain in $\lines$, there remain $\Theta(n)$ points in $\pts$. 
Thus, there are $\Theta(n^2)$ pairs of distinct points in $\pts^2$. 
The aforementioned process depends on whether there exists a point of $\pts$ that participates in $\Theta(n)$ horizontal pairs or not.

First, we consider the case where there exists a point $u\in \pts$ that participates in $\Theta(n)$ horizontal pairs. Let $\pts_u$ be the set of points of $\pts$ that form a horizontal pair with $u$.   
Since $D(\pts)=\Theta(\sqrt{n})$, the points of $\pts_u$ have $O(\sqrt{n})$ distinct $x$-coordinates.
Every $x$-coordinate repeats $O(\sqrt{n})$ times (otherwise, $D(\pts_u)$ would be too large).
Thus, there are $\Theta(\sqrt{n})$ subsets of points from $\pts_u$, such that all points in the same subset have the same $x$-coordinate.
Consider one such subset and let $\lines'$ be the set of $\Theta(\sqrt{n})$ lines of $\lines$ that participate in it. 
Since there is a point with the same $x$-coordinate on every line of $\lines'$, Theorem \ref{th:Freiman} implies that the $y$-intercepts in $\lines'$ form a generalized arithmetic progression of constant dimension and size $\Theta(\sqrt{n})$.

We now consider the case where no point of $\pts$ participates in $\Theta(n)$ horizontal pairs.
This implies that there are $o(n^2)$ horizontal pairs.
Since the total number of pairs is $\Theta(n^2)$, there are $\Theta(n^2)$ vertical pairs.
We say that two lines of $\lines$
are \emph{connected} if there is a vertical pair with one point on each line. 
Since any pair of lines of $\lines$ yields $O(n)$ vertical pairs, there are $\Theta(n)$ connected pairs of lines.
For any connected pair of lines, the difference of the two corresponding $y$-intercepts is a distance spanned by $\pts$.
Since every such difference of $y$-intercepts corresponds to  $O(\sqrt{n})$ pairs of connected lines, there are $\Theta(\sqrt{n})$ differences that repeat for $\Theta(\sqrt{n})$ pairs of connected lines. 

Let $Y$ be the set of $y$-intercepts of the lines of $\lines$.
Each difference that repeats $\Theta(\sqrt{n})$ times contributes $\Theta(n)$ quadruples to $E(Y)$. 
By the preceding paragraph, there are $\Theta(\sqrt{n})$ differences that repeat $\Theta(\sqrt{n})$ times. 
Thus, 
\[ E(Y)= \Theta(n^{3/2}) = \Theta(|Y|^3). \] 

By Theorem \ref{th:schoen}, there exists $Y'\subset Y$ of size $\Theta(\sqrt{n})$ such that $|Y'-Y'|=\Theta(|Y'|)$. 
Then, Theorem \ref{th:Freiman} implies that $Y'$ is contained in a generalized arithmetic progression of constant dimension and size $\Theta(\sqrt{n})$. 

To recap, in either case we find a set of $\Theta(\sqrt{n})$ lines of $\lines$ with $y$-intercepts that are contained in the a progression, as required.
After finding such a set, we remove it from $\lines$. We also remove the points that are on these lines from $\pts$. 
After $\Theta(1)$ steps of this process, we remain with $o(n)$ points of $\pts$ and the process ends.
We discard the remaining points.
\end{proof}

The following construction demonstrates that sets with few distances may be quite different than Cartesian products.
\vspace{2mm}

\noindent {\bf Theorem \ref{th:Construction}.}
\emph{There exists a set $\pts$ of $n$ points such that no two points of $\pts$ have the same $x$-coordinate and $D(\pts)=2\sqrt{n}-2$ under the $\ell_\infty$ metric.}

\begin{proof}
Let $\lines$ be the set of horizontal lines defined by $y=a$ with $a\in \{1,2,3,\ldots,\sqrt{n}\}$.
Let $B=\{b_1,b_2,b_3,\ldots,b_{\sqrt{n}}\}$ be a set of irrational numbers with the following properties. Every $b_i\in B$ satisfies $0<b_i<0.5$.
For every distinct $b_i,b_j\in B$, the difference $b_i-b_j$ is irrational. 

We construct $\pts$ by taking $\sqrt{n}$ points from every line of $\lines$, as follows. From the line defined by $y=a$, we add to $\pts$ the points $((b_a+a')/10n,a)$ with $a'\in \{1,2,3,\ldots,\sqrt{n}\}$.
In other words, we take an arithmetic progression starting at $(b_a+1)/10n$ and with step size $1/10n$.
See Figure \ref{fi:Moaaz} for an example. 

\begin{figure}[h]
    \centering
    \includegraphics[width=0.5\linewidth]{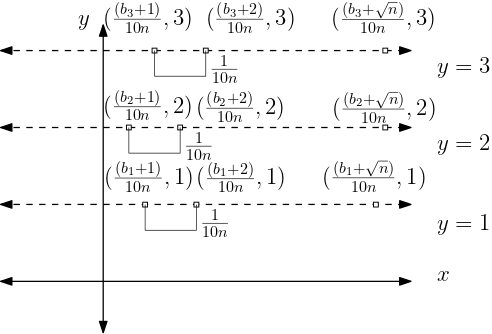}
    \caption{\small \sf Placing disjoint arithmetic progressions with step $1/10n$ on horizontal lines.}
\label{fi:Moaaz}
\end{figure}

Since we took $\sqrt{n}$ points from $\sqrt{n}$ lines, we obtain that $|\pts|=n$.
Since every two $b_i,b_j\in B$ satisfy that $b_i-b_j$ is irrational, no two numbers in $\pts$ have the same $x$-coordinate.
It remains to show that $D(\pts)=2\sqrt{n}-2$.
The distance between any two points on the same line of $\lines$ is in \[ \{1/10n,2/10n,3/10n,\ldots,(\sqrt{n}-1)/10n\}. \]
Thus, pairs of points from the same line contribute $\sqrt{n}-1$ distances.

Let $p,q\in \pts$ be points from different lines of $\lines$. Then $|p_y-q_y|\ge 1$ and $|p_x-q_x|<1$ (recall that all elements of $B$ are between 0 and 0.5).
By the definition of the $\ell_\infty$ metric, the distance between $p$ and $q$ is $|p_y-q_y|$.
By the definition of $\lines$, this distance is in $\{1,2,3,\ldots,\sqrt{n}-1\}$.
That is, points from different lines span $\sqrt{n}-1$ distances.

By considering points from the same line and points from different lines, we obtain $D(\pts)=2\sqrt{n}-2$. This concludes the proof.
\end{proof}

Finally, we provide a proof sketch for Corollary \ref{co:structureSubset}.
\vspace{2mm}

\noindent {\bf Corollary \ref{co:structureSubset}.}
\emph{Let $\pts$ be a set of $n$ points such that $D(\pts)=\Theta(\sqrt{n})$ under the $\ell_\infty$ metric.
Then there exist a subset $\pts' \subseteq \pts$ of $\Theta(n)$ points with the following properties:
\begin{itemize}[noitemsep,topsep=1pt,leftmargin=5mm]
    \item There exists a set $\lines$ of $\Theta(\sqrt{n})$ lines, either all vertical or all horizontal. Every line $\ell\in \lines$ satisfies $|\ell\cap \pts'| = \Theta(\sqrt{n})$. 
    \item There exists a generalized arithmetic progression $A$ with of dimension $\Theta(1)$ and size $\Theta(\sqrt{n})$. For every $\ell\in\lines$, when thinking of $\ell$ as $\RR$, there exists $r\in \RR$ such that $(\ell\cap \pts') \subseteq (r+A)$.
    \item The intercepts of the lines (with the orthogonal axis) are contained in another generalized arithmetic progression of dimension $\Theta(1)$ and size $\Theta(\sqrt{n})$.  
    \end{itemize}
}

\begin{proof}[Proof sketch.]
We apply Theorem \ref{th:Structure}(b), obtaining a set of lines $\lines$ and $\Theta(1)$ progressions $A_1,\ldots,A_s$. By the pigeonhole principle, there exists $A_j$ that corresponds to $\Theta(\sqrt{n})$ lines of $\lines$. That is, for each such line $\ell$ there exists $r\in \RR$ such that $|(r+A_j)\bigcap (\pts\cap \ell)| = \Theta(\sqrt{n})$.
We discard the other lines from $\lines$ and set $A=A_j$.
For each remaining line $\ell\in \lines$, we take an $r\in \RR$ such that $|(r+A)\bigcap (\pts\cap \ell)| = \Theta(\sqrt{n})$. We then add to $\pts'$ the points of $(r+A)\bigcap (\pts\cap \ell)$.

To recap, we remain with a set of $\Theta(\sqrt{n})$ lines, each containing $\Theta(\sqrt{n})$ points of $\pts'$. The points on each line are contained in a translation of $A$.
Note that $|\pts'|=\Theta(n)$. 

We repeat the process in the last part of the proof of Theorem \ref{th:Structure}(b).
This yields $\Theta(1)$ disjoint subsets of the lines of $\lines$, each of size $\Theta(\sqrt{n})$. In each subset, the intercepts of the lines (with the orthogonal axis) are contained in a generalized arithmetic progression of dimension $\Theta(1)$ and size $\Theta(\sqrt{n})$.
We arbitrarily choose one of these sets, discard the other lines from $\lines$, and remove the points of the discarded lines from $\pts'$. 
It is not difficult to verify that we still have $|\pts'| = \Theta(n)$.
\end{proof}

\parag{Acknowledgements.}
We are grateful to our mentors Adam Sheffer, N\'ora Frankl, Surya Mathialagan, and Jonathan Passant for their guidance and support and for making Polymath REU possible amidst the pandemic.


\begin{thebibliography}{99}
%
\bibitem{BELMPRT19}
V.\ Balaji, O.\ Edwards, A.\ M.\ Loftin, S.\ Mcharo, L.\ Phillips, A.\ Rice, and B.\ Tsegaye, Sets in Rd Determining k Taxicab Distances, \emph{Involve} {\bf 13} (2020), 487--509.
%
\bibitem{BKT04}
J.\ Bourgain, N.\ Katz, and T.\ Tao, A sum–product estimate in finite fields, and applications, 
\emph{Geom.\ Funct.\ Anal.} {\bf 14} (2004), 27--57.
%
\bibitem{CLOu}
D.\ Cox, J.\ Little, and D.\ O'Shea,
{\it Ideals, Varieties, and Algorithms: An Introduction to
Computational Algebraic Geometry and Commutative Algebra}, 3rd edition,
Springer-Verlag, Heidelberg, 2007.
%
\bibitem{erd46}
P.\ Erd\H os,
On sets of distances of $n$ points,
\emph{Amer. Math. Monthly} {\bf 53} (1946), 248--250.
%
\bibitem{Garibaldi04}
J.\ S.\ Garibaldi, 
Erd\H os distance problem for convex metrics, Ph.D. dissertation, University of California, Los Angeles, 2004.
%
\bibitem{GIS11}
J.\ Garibaldi, A.\ Iosevich, and S.\ Senger, 
\emph{The Erd\H os distance problem},
American Mathematical Soc., 2011.
%
\bibitem{GK15}	
L.\ Guth and N.H.\ Katz,
On the Erd{\H o}s distinct distances problem in the plane,
{\em Annals Math.} {\bf 181} (2015), 155--190.
%
\bibitem{IL04}
A.\ Iosevich and I.\ \L{}aba. Distance sets of well-distributed planar point sets, 
\emph{Discrete \& Computational Geometry} {\bf 31} (2004), 243--250.
%
\bibitem{Mat11}
J.\ Matou\v sek, 
The number of unit distances is almost linear for most norms,
\emph{Advances in Mathematics} {\bf 226} (2011), 2618--2628.
%
\bibitem{MPPRS03}
B.\ Murphy, G.\ Petridis, T.\ Pham, M.\ Rudnev, and S.\ Stevens, On the Pinned
Distances Problem over Finite Fields, arXiv:2003.00510.
%
\bibitem{Schoen15}
T.\ Schoen, 
New bounds in Balog--Szemer\'edi--Gowers theorem,
\emph{Combinatorica} {\bf 35} (2015), 695--701.
%
\bibitem{SZ16}
M.\ Sharir and J.\ Zahl,
Cutting algebraic curves into pseudo-segments and applications,
\emph{J.\ Combinat.\ Theory Ser.\ A} {\bf 150} (2017), 1--35.
%
\bibitem{Sheffer14}
A.\ Sheffer, Distinct Distances: Open Problems and Current Bounds,
arXiv:1406.1949.
%
\bibitem{SheZa20}
A.\ Sheffer and J.\ Zahl,
Distinct distances in the complex plane, arXiv:2006.08886.
%
\bibitem{SoTo01}
J.\ Solymosi and C.\ D.\ T\'oth, 
Distinct distances in the plane, 
\emph{Discrete Comput. Geom.} {\bf 25} (2001), 629--634.
%
\bibitem{Szek97}
L.\ Sz\'ekely, Crossing numbers and hard Erd\H os problems in discrete geometry, \emph{Combinatorics,
Probability and Computing} {\bf 6} (1997), 353--358.
%
\bibitem{TV06}
T.\ Tao and V.\ H.\ Vu,
\emph{Additive combinatorics}, Cambridge University Press, 2006.
%
\bibitem{Tardos03}
G.\ Tardos, 
On distinct sums and distinct distances, 
\emph{Advances Math.} {\bf 180} (2003), 275--289.
%
\end{thebibliography}
\end{document}